\documentclass[12pt,sort&compress]{elsarticle}

\usepackage{bm}

\usepackage{mathrsfs}
\usepackage{amsmath}
\usepackage{amsfonts}
\usepackage{amsthm}
\usepackage{color}
\usepackage{xcolor}
\usepackage{caption}
\usepackage{subcaption}
\usepackage{float}
\usepackage{url}
\usepackage{multicol}
\usepackage[top=1in,bottom=1in,left=1in,right=1in]{geometry}
\usepackage[pdfborder={0 0 0},colorlinks,allcolors=blue]{hyperref}
\usepackage{txfonts}
\usepackage{setspace}
\usepackage{arydshln}
\usepackage{enumitem}
\usepackage{lipsum}% for dummy text
\usepackage{siunitx,booktabs}
\usepackage{nth}
\usepackage{accents} % Accents: circle

\theoremstyle{definition}%plain/definition/remark

\newtheorem{theorem}{Theorem}%[section]
\newtheorem{lemma}{Lemma}%[section]
\newtheorem{assumption}{Assumption}%[section]
\theoremstyle{definition}%plain/definition/remark
\newtheorem{remark}{Remark}%[section]

\allowdisplaybreaks

\newcommand{\vertiii}[1]{{\left\vert\kern-0.2ex\left\vert\kern-0.2ex\left\vert #1 
    \right\vert\kern-0.2ex\right\vert\kern-0.2ex\right\vert}}

\newcommand{\vertii}[1]{{\left\vert\kern-0.2ex\left\vert #1 
    \right\vert\kern-0.2ex\right\vert}}

\usepackage{listings}
\definecolor{light-gray}{gray}{0.95}
\begin{document}
\begin{frontmatter}

\title{Variational multiscale modeling with discretely divergence-free subscales: Non-divergence-conforming discretizations}

\author[ucsd]{Sajje~Lee~Calfy}
\author[colorado]{John~A.~Evans}
\author[ucsd]{David~Kamensky\corref{cor1}}
\ead{dmkamensky@eng.ucsd.edu}
\cortext[cor1]{Corresponding author}

\address[ucsd]{Department of Mechanical and Aerospace Engineering, University of California, San Diego, 9500 Gilman Drive, Mail Code 0411, La Jolla, CA 92093, USA}
\address[colorado]{Department of Aerospace Engineering Sciences, University of Colorado at Boulder, 429 UCB, Boulder, CO 80309, USA}

\journal{}

\begin{abstract}
A recent paper [J. A. Evans, D. Kamensky, Y. Bazilevs, ``Variational multiscale modeling with discretely divergence-free subscales'', Computers \& Mathematics with Applications, 80 (2020) 2517--2537] introduced a novel stabilized finite element method for the incompressible Navier--Stokes equations, which combined residual-based stabilization of advection, energetic stability, and satisfaction of a discrete incompressibility condition.  However, the convergence analysis and numerical tests of the cited work were subject to the restrictive assumption of a divergence-conforming choice of velocity and pressure spaces, where the pressure space must contain the divergence of every velocity function.  The present work extends the convergence analysis to arbitrary inf-sup-stable velocity--pressure pairs (while maintaining robustness in the advection-dominated regime) and demonstrates the convergence of the method numerically, using both the traditional and isogeometric Taylor--Hood elements.
\end{abstract}


\end{frontmatter}

\section{Introduction}\label{sec:introduction}

Our point of departure in the present work is \cite{Evans2019a}, which proposed a new stabilized finite element formulation for the incompressible Navier--Stokes equations.  The novelty of this formulation was its ability to combine residual-based stabilization of advection with energetic stability and satisfaction of a discrete incompressibility condition for a discrete velocity solution $\mathbf{u}^h$, i.e.,
\begin{equation}\label{eq:disc-div-free}
    \int_\Omega q^h\nabla\cdot\mathbf{u}^h\,d\Omega = 0\quad\forall q^h\in\mathcal{Q}^h\text{ ,}
\end{equation}
where $\Omega$ is the problem domain and $\mathcal{Q}^h$ is the discrete pressure space.  The basic Galerkin discretization of incompressible flow obviously satisfies \eqref{eq:disc-div-free}, but is unstable for the high Reynolds numbers ubiquitous in applications.  Widely-used residual-based stabilized methods such as Galerkin least squares (GLS) \cite{Franca1993} or the combination of streamline/upwind-Petrov--Galerkin (SUPG) \cite{brooks1982streamline} and pressure-stabilizing Petrov--Galerkin (PSPG) \cite{Hughes1986} typically fail to satisfy \eqref{eq:disc-div-free}, due to an inner product of the strong problem's momentum-balance residual with the gradient of the pressure test function.  The method of \cite{Evans2019a} circumvents this difficulty through a novel application of the variational multiscale (VMS) \cite{Hug98,HugSan06,Hughes07b,Bazilevs07b,Bazilevs09d} concept: a second discrete pressure field is introduced, formally taking the place of the fine-scale pressure in a VMS scale separation. 

The analysis of \cite{Evans2019a} showed that the proposed method satisfied energetic stability \cite[Lemma 2]{Evans2019a}, but only proved convergence of a simplified method for the Oseen equations---a linearized model of Navier--Stokes flow---under the restrictive assumption of a divergence-conforming discretization, i.e.,
\begin{equation}
\nabla\cdot\mathbf{v}^h\in\mathcal{Q}^h\quad\forall\mathbf{v}^h\in\mathcal{V}^h\text{ ,}
\end{equation}
where $\mathcal{V}^h$ is the discrete velocity space.  However, it was conjectured that the Oseen flow convergence result in fact holds for arbitrary choices of $\mathcal{V}^h \subset (H^1(\Omega))^d$ and $\mathcal{Q}^h \subset H^1(\Omega)$ satisfying the classic inf-sup stability condition for saddle point problems \cite{Brezzi1974,Babuska71,Boffi2008} (e.g., the Taylor--Hood element \cite{Taylor1973,Brezzi1991}, its isogeometric counterpart \cite{Hosseini2015}, the MINI element \cite{Arnold1984}, and the isogeometric subgrid method \cite{Bressan2013}).  The present paper shows that this conjecture is true, generalizing the convergence proof of \cite{Evans2019a} to non-divergence-conforming discretizations while maintaining robustness in the advection-dominated regime.  

We begin, in Section \ref{sec:formulation}, by restating the formulation proposed by \cite{Evans2019a}.  Section \ref{sec:oseen-conv} then gives the generalized convergence proof.  This convergence theorem is confirmed for the Oseen equations and empirically shown to extend to the full nonlinear Navier--Stokes equations in the numerical results of Section \ref{sec:numerical}.  We then draw conclusions and discuss potential future research directions in Section \ref{sec:conclusions}.

\section{Formulation}\label{sec:formulation}

In this section, we recall the VMS formulation proposed in \cite{Evans2019a}.  This formulation approximates the incompressible Navier--Stokes equations on a domain $\Omega\subset\mathbb{R}^d~,~d\in\{2,3\}$, with a Newtonian viscosity model.  For simplicity, we state the problem and analyze our formulation assuming homogeneous Dirichlet boundary conditions on velocity.  In strong form, this partial differential equation (PDE) system is:  Find velocity $\mathbf{u}:\overline{\Omega}\times\lbrack 0,T\rbrack\to\mathbb{R}^d$ and pressure $p:\overline{\Omega}\times(0,T)\to\mathbb{R}$ such that
\begin{equation}
    \begin{array}{rclr}\partial_t\mathbf{u} + \mathbf{u}\cdot\nabla\mathbf{u} - \nabla\cdot\left(2\nu\nabla^s\mathbf{u}\right) + \nabla p &=& \mathbf{f} & \text{in }\Omega\times (0,T)\\
    \nabla\cdot\mathbf{u} &=& 0 & \text{in }\Omega\times(0,T)\\
    \mathbf{u}&=&\bm{0}&\text{on }\partial\Omega\times(0,T)\\
    \mathbf{u}\vert_{t=0} &=& \mathbf{u}_0&\text{in }\Omega\text{ ,}\end{array}\label{eq:strong-problem}
\end{equation}
where $\nu > 0$ is the kinematic viscosity, $\nabla^s$ is the symmetrized gradient operator, $\mathbf{f}$ is the body force per unit volume, and $T > 0$ is the time over which the system evolves, beginning from the initial velocity field $\mathbf{u}_0$.  

Finite element methods are based on a weakened, variational form of this problem.  To state variational problems succintly, we introduce some notational conventions.  First, inner products of functions written with no subscript are to be understood as $L^2(\Omega)$ inner products, e.g.,
\begin{equation}
(u,v) = (u,v)_{L^2(\Omega)} = \int_\Omega uv\,d\Omega
\end{equation}
and likewise for the corresponding induced norm, $\Vert u\Vert = \sqrt{(u,u)}$.  Further, in discrete problems where $\Omega$ is partitioned into a mesh of elements $\{\Omega_e\}_{e=1}^{n_{\text{el}}}$ with disjoint interiors and closures covering the domain, the $L^2$ inner product is to be interpreted in a ``broken'' sum-over-elements sense, excluding any singular distributional parts of its arguments at interior boundaries between elements, e.g.,
\begin{equation}
    (u,v) = \sum_{e=1}^{n_{\text{el}}}(u,v)_{L^2(\Omega_e)}
\end{equation}
(and, again, likewise for the norm), which of course maintains backward compatibility with the usual definition when functions are in $L^2(\Omega)$.  With this notation in hand, our weak problem is:  Find $(\mathbf{u},p)\in\mathcal{X}_T$ with $\mathbf{u}(0)=\mathbf{u}_0$ such that, $\forall(\mathbf{v},q)\in\mathcal{X}$ and at a.e. $t\in (0,T)$,
\begin{equation}
    \left(\partial_t\mathbf{u}(t),\mathbf{v}\right) + c(\mathbf{u}(t),\mathbf{u}(t),\mathbf{v}) + k(\mathbf{u}(t),\mathbf{v}) - b(\mathbf{v},p(t)) + b(\mathbf{u}(t),q) = (\mathbf{f}(t),\mathbf{v})\text{ ,}\label{eq:ns-weak}
\end{equation}
where the velocity, pressure, and mixed function spaces are defined as
\begin{equation}
    \mathcal{V}:=\left(H_0^1(\Omega)\right)^d\quad\text{,}\quad        \mathcal{Q}:=L^2_0(\Omega)\quad\text{,}\quad
        \mathcal{X}:= \mathcal{V}\times\mathcal{Q}\text{ ,}
\end{equation}
with time-dependent counterparts
\begin{align}
    \mathcal{V}_T&:=\left\{\mathbf{v}\in C^0\left(\lbrack 0,T\rbrack;\mathcal{V}\right)~:~\partial_t\mathbf{v}\in L^2((0,T);\mathcal{V})\right\}\text{ ,}\\
        \mathcal{Q}_T&:=L^2((0,T);\Omega)\text{ ,}\\
        \mathcal{X}_T&:= \mathcal{V}_T\times\mathcal{Q}_T\text{ ,}
\end{align}
and the convective, diffusive, and constraint forms are defined by
\begin{align}
    c(\mathbf{v}_1,\mathbf{v}_2,\mathbf{v}_3) &= \int_{\Omega}(\mathbf{v}_1\cdot\nabla\mathbf{v}_2)\cdot\mathbf{v}_3\,d\Omega\text{ ,} \nonumber \\
    k(\mathbf{v}_1,\mathbf{v}_2) &= \int_\Omega 2\nu\nabla^s\mathbf{v}_1 : \nabla^s\mathbf{v}_2\,d\Omega\text{ ,} \nonumber \\
    b(\mathbf{v},q) &= \int_\Omega\nabla\cdot\mathbf{v} q\,d\Omega\text{ .} \nonumber
\end{align}
To state our numerical method, we also introduce two alternative convective forms, which are equivalent to the $c$ defined above when its first argument is exactly solenoidal:
\begin{align}
    c_{\text{cons}}(\mathbf{v}_1,\mathbf{v}_2,\mathbf{v}_3) &= -\int_\Omega\mathbf{v}_2\cdot\left(\mathbf{v}_1\cdot\nabla\mathbf{v}_3\right)\,d\Omega\text{ ,}\\
    c_{\text{skew}}(\mathbf{v}_1,\mathbf{v}_2,\mathbf{v}_3) &= \frac{1}{2}\left(c(\mathbf{v}_1,\mathbf{v}_2,\mathbf{v}_3) + c_{\text{cons}}(\mathbf{v}_1,\mathbf{v}_2,\mathbf{v}_3)\right)\text{ .}
\end{align}
\noindent
The semi-discrete formulation of the method is developed in \cite[Section 2.5]{Evans2019a}, based on VMS reasoning, and uses the following function spaces:
\begin{align}
    \mathcal{V}^h_T&:=\left\{\mathbf{v}\in C^0(\lbrack 0,T\rbrack;\mathcal{V}^h)~:~\partial_t\mathbf{v}\in L^2((0,T);\mathcal{V}^h)\right\}\text{ ,}\\
    \mathcal{V}'_T&:=\left\{\mathbf{v}\in C^0(\lbrack 0,T\rbrack;\mathcal{V}')~:~\partial_t\mathbf{v}\in L^2((0,T);\mathcal{V}')\right\}\text{ ,}\\
    \mathcal{Q}^h_T&:= L^2((0,T);\mathcal{Q}^h)\text{ ,}\\
    \mathcal{Q}'_T&:= L^2((0,T);\mathcal{Q}')\text{ ,}
\end{align}
\noindent
where $\mathcal{V}^h$ and $\mathcal{Q}^h$ are finite-dimensional subsets of $\mathcal{V}$ and $\mathcal{Q}\cap H^1(\Omega)$, which are assumed to satisfy an appropriate inf-sup condition, and $\mathcal{V}'$ and $\mathcal{Q}'$ are the fine-scale velocity and pressure spaces.  In standard piecewise-polynomial finite element spaces, the assumption that $\mathcal{Q}^h\subset H^1(\Omega)$ corresponds to a restriction that the discrete pressure space must be $C^0$.  We assume that $\mathcal{Q}' = \mathcal{Q}^h$,\footnote{We note that this is at odds with the VMS direct-sum decomposition of spaces into coarse and fine components.  However, this work follows the philosophy of using VMS as a formal tool to inspire new stabilized methods, which must then be proven correct through further analysis.} and will formally eliminate $\mathcal{V}'$ through static condensation.  The semi-discrete formulation is then:  Find $(\mathbf{u}^h,p^h)\in\mathcal{V}_T^h\times\mathcal{Q}_T^h$ and $(\mathbf{u}',p')\in\mathcal{V}_T'\times\mathcal{Q}_T'$ satisfying $\mathbf{u}^h(0)=\mathbf{u}^h_0$, $\mathbf{u}'(0)=\mathbf{u}'_0$, such that $\forall (\mathbf{v}^h,q^h)\in\mathcal{X}^h$ and a.e. $t\in (0,T)$,
\begin{align}
    \nonumber &(\partial_t\mathbf{u}^h(t),\mathbf{v}^h) + c_{\text{skew}}(\mathbf{u}^h(t),\mathbf{u}^h(t),\mathbf{v}^h) + k(\mathbf{u}^h(t),\mathbf{v}^h) - b(\mathbf{v}^h,p^h(t)) + b(\mathbf{u}^h(t),q^h)\\
    \nonumber &+c_{\text{cons}}(\mathbf{u}^h(t),\mathbf{u}'(t),\mathbf{v}^h)+c_{\text{skew}}(\mathbf{u}'(t),\mathbf{u}^h(t),\mathbf{v}^h)+c_{\text{cons}}(\mathbf{u}'(t),\mathbf{u}'(t),\mathbf{v}^h)+(\partial_t\mathbf{u}'(t),\mathbf{v}^h)\\
    &+ (\tau_\text{C}\nabla\cdot\mathbf{u}^h(t),\nabla\cdot\mathbf{v}^h) =(\mathbf{f}(t),\mathbf{v}^h) \label{eq:semidisc-coarse}
\end{align}
and $\forall (\mathbf{v}',q')\in\mathcal{V}'\times\mathcal{Q}'$ and a.e. $t\in (0,T)$,
\begin{align}
    \nonumber & (\partial_t\mathbf{u}'(t),\mathbf{v}') + (\tau^{-1}_{\text{M}}\mathbf{u}'(t),\mathbf{v}') + c(\mathbf{u}'(t),\mathbf{u}^h(t),\mathbf{v}')\\
    &+(\nabla p'(t),\mathbf{v}') - (\nabla q',\mathbf{u}'(t)) = -(\mathbf{r}_\text{M}(t),\mathbf{v}') \text{ ,} \label{eq:semidisc-fine}
\end{align}
 where $\tau_\text{M}$ and $\tau_\text{C}$ are stabilization parameters.  For the purposes of this paper's analysis, we assume the definitions
\begin{equation}\label{eq:stab-param-defn-ns}
\tau_\text{M}=\min{\left\{\frac{h}{2\left|\textbf{u}^h\right|},\frac{h^2}{C_{\text{inv}}\nu}\right\}}\quad\text{and}\quad\tau_\text{C} = \max\{h\vert\mathbf{u}^h\vert,\nu\}\text{ ,}
\end{equation}
on a quasi-uniform mesh with global element size $h$, where $C_{\text{inv}}>0$ is a constant independent of $h$, $\nu$, and $\mathbf{u}^h$, related to inverse estimates (cf. Assumption \ref{assumption2} in the sequel).  More complicated definitions for nonuniform and/or anisotropic meshes are available in the literature, with the same general asymptotic behavior in the advective and diffusive limits, e.g., \cite[(63)--(69)]{Bazilevs07b}.\footnote{Some authors consider alternative asymptotic behavior of $\tau_\text{M}$, e.g., \cite[(3.2)]{Matthies2009} (using a different notation), but we restrict our analysis to the classical choices in \eqref{eq:stab-param-defn}.  See also, Remark \ref{rem:different-tau_M-asymptotics} in the sequel.}

\begin{remark}
The semidiscrete formulation \eqref{eq:semidisc-coarse}--\eqref{eq:semidisc-fine} was shown to be energetically-stable for the full Navier--Stokes problem in \cite[Lemma 2]{Evans2019a}.  The careful selection of advective forms $c$, $c_\text{cons}$ and $c_\text{skew}$ in different terms was crucial to obtaining this result.
\end{remark}

We discretize this problem in time using an implicit midpoint rule.  In particular, this allows the fine-scale velocity subproblem to be statically condensed, i.e., the fine-scale velocity at time step $n$ is solved for pointwise, to obtain
\begin{equation}
\textbf{u}_n^\prime=\left(\left(\frac{1}{\Delta t}+\frac{1}{2\tau_\text{M}}\right)\mathbf{I}+\frac{1}{2}\nabla\textbf{u}_{n-\frac{1}{2}}^h\right)^{-1}\left(-\textbf{r}_M\left(t_{n-1/2}\right)-\nabla p^\prime+\left(\left(\frac{1}{\Delta t}-\frac{1}{2\tau_\text{M}}\right)\mathbf{I}-\frac{1}{2}\nabla\textbf{u}_{n-\frac{1}{2}}^h\right)\textbf{u}_{n-1}^\prime\right)\text{ ,}\label{eq:static_cond}
\end{equation}
where
\begin{align}
\textbf{u}_{n-\frac{1}{2}}^h = \frac{1}{2}\left(\textbf{u}_n^h+\textbf{u}_{n-1}^h\right)\quad\text{and}\quad
\textbf{u}_{n-\frac{1}{2}}^\prime = \frac{1}{2}\left(\textbf{u}_n^\prime+\textbf{u}_{n-1}^\prime\right)
\end{align}
are the midpoint coarse- and fine-scale velocities and $\Delta t$ is the time step size.  We also explore a simplified model, in which a quasi-static approximation is used for the fine-scale velocity:
\begin{equation}
    \mathbf{u}' = -\tau_{\text{M}}\left(\nabla p' + \mathbf{r}_\text{M}\right)\quad\text{with}\quad\tau_\text{M} = \min{\left\{\frac{\Delta t}{2},\frac{h}{2\left|\textbf{u}^h\right|},\frac{h^2}{C_{\text{inv}}\nu}\right\}}\text{ .}
\end{equation}
We find that, in practice, this simplified model performs well in most situations, but it is not covered by the energy stability theory of \cite{Evans2019a} and may, in principle, exhibit stability problems when the time step is extremely small relative to the spatial element diameter $h$, i.e., $\vert\mathbf{u}\vert\Delta t \ll h$.  The small time step issue is discussed further in \cite{Hsu2010} and references cited therein.

\section{Proof of convergence for the Oseen problem}\label{sec:oseen-conv}

In this section we show convergence of our method applied to the model problem of steady Oseen flow, in which we neglect time dependence and replace the advection velocity (i.e., first argument to $c$ in \eqref{eq:ns-weak}) with a known solenoidal vector field, resulting in a linear problem. The weak form of the Oseen problem is: Find $\textbf{u}\in\mathcal{V}$ and $p\in\mathcal{Q}$ such that for all $\textbf{v}\in\mathcal{V}$ and $q\in\mathcal{Q}$,
\begin{equation}
c\left(\textbf{a},\textbf{u},\textbf{v}\right)+k\left(\textbf{u},\textbf{v}\right)-b\left(\textbf{v},p\right)+b\left(\textbf{u},q\right)=\left(\textbf{f},\textbf{v}\right)\text{ ,}
\end{equation}
where the advection velocity $\textbf{a}$ satisfies $\nabla\cdot\textbf{a}=0$. For simplicity of analysis, we assume that $\mathbf{a}$ is a constant vector throughout the domain.  Analysis based on this assumption can be extended to variable $\mathbf{a}$ through additional bookkeeping, as in, e.g., \cite{Franca1993}.  Thus, a discretization for our simplified model problem analogous to the semi-discrete problem \eqref{eq:semidisc-coarse}--\eqref{eq:semidisc-fine} is: Find $(\mathbf{u}^h,p^h,\mathbf{u}',p')\in\mathcal{V}^h\times\mathcal{Q}^h\times\mathcal{V}'\times\mathcal{Q}'$, such that for all $(\mathbf{v}^h,q^h,\mathbf{v}',q')\in\mathcal{V}^h\times\mathcal{Q}^h\times\mathcal{V}'\times\mathcal{Q}'$,
\begin{equation}
    A\left((\mathbf{u}^h,p^h,\mathbf{u}',p'),(\mathbf{v}^h,q^h,\mathbf{v}',q')\right) = (\mathbf{f},\mathbf{v}^h+\mathbf{v}')\text{ ,}
\end{equation}
where:
\begin{align}
    \nonumber & A\left((\mathbf{u}^h,p^h,\mathbf{u}',p'),(\mathbf{v}^h,q^h,\mathbf{v}',q')\right) \\
    \nonumber &~~~ =c(\mathbf{a},\mathbf{u}^h,\mathbf{v}^h) + k(\mathbf{u}^h,\mathbf{v}^h) - b(\mathbf{v}^h,p^h) + b(\mathbf{u}^h,q^h)\\
    \nonumber &~~~ +c_\text{cons}(\mathbf{a},\mathbf{u}',\mathbf{v}^h) + (\tau_\text{C}\nabla\cdot\mathbf{u}^h,\nabla\cdot\mathbf{v}^h) + (\tau_\text{M}^{-1}\mathbf{u}',\mathbf{v}')\\
    \nonumber&~~~-b(\mathbf{v}',p') + b(\mathbf{u}',q') + (\mathbf{a}\cdot\nabla\mathbf{u}^h,v')\\
    &~~~-(\nabla\cdot(2\nu\nabla^s\mathbf{u}^h),\mathbf{v}') + (\nabla p^h,\mathbf{v}')
\end{align}
and we replace $\mathbf{u}^h$ with $\mathbf{a}$ in the definitions of $\tau_\text{M}$ and $\tau_\text{C}$ given by \eqref{eq:stab-param-defn-ns}:
\begin{equation}\label{eq:stab-param-defn}
\tau_\text{M}=\min{\left\{\frac{h}{2\left|\textbf{a}\right|},\frac{h^2}{C_{\text{inv}}\nu}\right\}}\quad\text{and}\quad\tau_\text{C} = \max\{h\vert\mathbf{a}\vert,\nu\}\text{ .}
\end{equation}
Recalling that $\mathcal{Q}^h=\mathcal{Q}^\prime$, we introduce the notation
\begin{equation}
\tilde{\mathcal{Q}}=\mathcal{Q}^h+\mathcal{Q}^\prime = \mathcal{Q}^h 
\end{equation}
and define
\begin{equation}
\tilde{p}=p^h+p^\prime\quad\text{and}\quad\tilde{q}=q^h+q^\prime\quad\text{in}\quad\tilde{\mathcal{Q}}\text{ .}
\end{equation}
Finally, due to the absence of time derivatives, the fine-scale velocity static condensation \eqref{eq:static_cond} simplifies to
\begin{equation}
    \mathbf{u}' = \tau_\text{M}\left(\mathbf{f} - \mathbf{a}\cdot\nabla\mathbf{u}^h + \nabla\cdot(2\nu\nabla^s\mathbf{u}^h) - \nabla \tilde{p} \right)\text{ ,}
\end{equation}
and we obtain the following reduced formulation for the Oseen problem: Find $\left(\textbf{u}^h,p^h,\tilde{p}\right)\in\mathcal{V}^h\times\mathcal{Q}^h\times\tilde{\mathcal{Q}}$ such that for all $\left(\textbf{v}^h,q^h,\tilde{q}\right)\in\mathcal{V}^h\times\mathcal{Q}^h\times\tilde{\mathcal{Q}}$,
\begin{equation}
A_\text{red}\left(\left(\textbf{u}^h,p^h,\tilde{p}\right),\left(\textbf{v}^h,q^h,\tilde{q}\right)\right)=\left(\textbf{f},\textbf{v}^h + \tau_\text{M}\left(\mathbf{a}\cdot\nabla\mathbf{v}^h + \nabla\tilde{q}\right)\right)\text{ ,}
\end{equation}
where $A_\text{red}$ is given by
\begin{equation}
\begin{split}
A_\text{red}\left(\left(\textbf{u}^h,p^h,\tilde{p}\right),\left(\textbf{v}^h,q^h,\tilde{q}\right)\right) &= c\left(\textbf{a},\textbf{u}^h,\textbf{v}^h\right)+k\left(\textbf{u}^h,\textbf{v}^h\right)-b\left(\textbf{v}^h,p^h\right)+b\left(\textbf{u}^h,q^h\right) \\ 
&+\left(\textbf{a}\cdot\nabla\textbf{u}^h-\nabla\cdot\left(2\nu\nabla^s\textbf{u}^h\right)+\nabla\tilde{p},\tau_\text{M}\left(\textbf{a}\cdot\nabla\textbf{v}^h+\nabla\tilde{q}\right)\right) \\
&+\left(\tau_\text{C}\nabla\cdot\textbf{u}^h,\nabla\cdot\textbf{v}^h\right)\text{ .}
\end{split}\label{eq:A_red_defn}
\end{equation}
\noindent
Our analysis of this reduced problem will use the norms
\begin{equation}
\begin{split}
{\vertiii{(\textbf{u}^h,p^h,\tilde{p})}}^2 &= k\left(\textbf{u}^h,\textbf{u}^h\right)+\alpha\left(p^h,\ p^h\right)+\left(\tau_\text{C}\nabla\cdot\textbf{u}^h,\nabla\cdot\textbf{u}^h\right) \\
&+\left(\tau_\text{M}\left(\textbf{a}\cdot\nabla\textbf{u}^h+\nabla\tilde{p}\right),\ \textbf{a}\cdot\nabla\textbf{u}^h+\nabla\tilde{p}\right) \\
&+\left(\tau_\text{M}\nabla\cdot\left(2\nu\nabla^s\textbf{u}^h\right),\nabla\cdot\left(2\nu\nabla^s\textbf{u}^h\right)\right)
\end{split} \label{triple-bar-norm}
\end{equation}
and
\begin{equation}
{\vertiii{\left(\textbf{u}^h,p^h,\tilde{p}\right)}}_\ast^2={\vertiii{\left(\textbf{u}^h,p^h,\tilde{p}\right)}}^2+{\left\Vert\tau_\text{M}^{-\frac{1}{2}}\textbf{u}^h\right\Vert}^2+{\left\Vert\tau_\text{C}^{-\frac{1}{2}}p^h\right\Vert}^2\text{ ,} \label{triple-bar-norm-star}
\end{equation}
where
\begin{equation}
\alpha=\nu^{-1}\min{\left\{1,{\rm Pe}^{-2},{\rm Pe}_h^{-1}\right\}}\text{ ,} \label{Define:alpha}
\end{equation}
\noindent
in which ${\rm Pe}$ and ${\rm Pe}_h$ are the global and element P\'{e}clet numbers, 
\begin{equation}
{\rm Pe}=\frac{\vert\textbf{a}\vert L}{2\nu}\quad\text{and}\quad {\rm Pe}_h=\frac{\vert\textbf{a}\vert h}{2\nu}\text{ ,} \label{Define:PecletNumbers}
\end{equation}
and $L$ is a global length scale associated with domain $\Omega$.

Before proceeding with our convergence proof, we make the following assumptions of the discrete function spaces:
\begin{assumption}\label{assumption1}
For any $q^h\in\mathcal{Q}^h$, there exists $\textbf{w}^h\in\mathcal{V}^h$ such that:
\begin{equation}
-b\left(\textbf{w}^h,q^h\right)\geq\tilde{\gamma}{\vertii{\nabla\textbf{w}^h}}\,{\vertii{q^h}}\text{ ,}
\end{equation}
where $\tilde{\gamma} > 0$ is an $h$-independent constant and $\textbf{w}^h$ can be scaled such that
\begin{equation}\label{eq:norm-grad-w}
{\vertii{\nabla\textbf{w}^h}}={\vertii{q^h}}\text{ .}
\end{equation}
This is true whenever $\mathcal{V}^h$ and $\mathcal{Q}^h$ satisfy a standard inf-sup condition for incompressible flow.
\end{assumption}

Standard piecewise polynomial finite element and isogeometric function spaces also typically satisfy the following inverse inequality \cite{Bazilevs2006}.
\begin{assumption}\label{assumption2}
For each $\Omega^e$, there exists a constant $C_\text{inv} > 0$ independent of the diameter $h_e$ of $\Omega^e$ such that:
\begin{equation}
\left(\nabla\cdot\left(2\nu\nabla^s\textbf{v}^h\right),\nabla\cdot\left(2\nu\nabla^s\textbf{v}^h\right)\right)\le\frac{C_\text{inv}}{h_e^2}\left(\nu\nabla^s\textbf{v}^h,\nu\nabla^s\textbf{v}^h\right)\label{eq:inv-est}
\end{equation}
for all $\textbf{v}^h\in\mathcal{V}^h$.  On a quasi-uniform mesh, $C_\text{inv}$ can be considered a global constant.
\end{assumption}

We begin our convergence proof by noting that, due to the residual-based nature of the stabilization, we have the following consistency result for the bilinear form $A_\text{red}$:
\begin{lemma}[Consistency]\label{Consistency}
If the exact solution $(\mathbf{u},p)\in \left(\mathcal{V}\cap(H^2(\Omega))^d\right)\times \left(\mathcal{Q}\cap H^1(\Omega)\right)$, then the reduced formulation is strongly consistent:
\begin{equation}
A_\text{red}\left(\left(\textbf{u}-\textbf{u}^h,p-p^h,p-\tilde{p}\right),\left(\textbf{v}^h,q^h,\tilde{q}\right)\right)=0~~~~~~~ \forall\left(\textbf{v}^h,q^h,\tilde{q}\right)\in\mathcal{V}^h\times\mathcal{Q}^h\times\tilde{\mathcal{Q}}\text{ .}
\end{equation}
\end{lemma}
We next prove a useful bound on our bilinear form. 
\begin{lemma}[Continuity]\label{Continuity}
The bilinear form $A_\text{red}$ satisfies
\begin{equation}
A_\text{red}\left(\left(\textbf{u},p,\tilde{p}\right),\left(\textbf{v}^h,q^h,\tilde{q}\right)\right)\le 6{\vertiii{\left(\textbf{u},p,\tilde{p}\right)}}_\ast{\vertiii{\left(\textbf{v}^h,q^h,\tilde{q}\right)}}\label{eq:A-bdd}
\end{equation}
for all $\left(\textbf{u},p,\tilde{p}\right)\in{\mathring{\mathcal{V}}_{\mathcal{Q}^h}\cap (H^2(\Omega)^d}+\mathcal{V}^h)\times\mathcal{Q}\cap H^1(\Omega)\times\mathcal{Q}\cap H^1(\Omega)$ and $\left(\textbf{v}^h,q^h,\tilde{q}\right)\in\mathcal{V}^h\times\mathcal{Q}^h\times\tilde{\mathcal{Q}}$, where
\begin{equation}
    \mathring{\mathcal{V}}_{\mathcal{Q}^h} := \left\{\mathbf{w}\in\mathcal{V}~:~b\left(\mathbf{w},q^h\right) = 0~\forall q^h\in\mathcal{Q}^h\right\}\text{ ,}
\end{equation}
i.e., the subspace of $\mathcal{V}$ that is discretely divergence-free with respect to $\mathcal{Q}^h$.
\begin{proof}
To prove the desired result, we bound each term of $A_\text{red}\left(\left(\textbf{u},p,\tilde{p}\right),\left(\textbf{v}^h,q^h,\tilde{q}\right)\right)$, as expanded in the right-hand side of \eqref{eq:A_red_defn}.\footnote{N.b. that $(\mathbf{u},p,\tilde{p})$ here is quantified over the space stated in the lemma, and does not refer to the exact solution.}

\noindent
$\textbf{Continuity of Term 1}$: The first term of $A_\text{red}\left(\left(\textbf{u},p,\ \tilde{p}\right),\left(\textbf{v}^h,q^h,\tilde{q}\right)\right)$ is
\begin{equation}
c\left(\textbf{a},\textbf{u},\textbf{v}^h\right) = \int_{\Omega}{\left(\textbf{a}\cdot\nabla\textbf{u}\right)\cdot\textbf{v}^h\,d\Omega}\text{ .}\label{eq:term-1}
\end{equation}
\noindent
Since $\mathcal{Q}^h=\mathcal{Q}^\prime$ and $\mathbf{u}\in\mathring{\mathcal{V}}_{\mathcal{Q}^h}$,
\begin{equation}
c\left(\textbf{a},\textbf{u},\textbf{v}^h\right)=\int_{\Omega}{\left(\textbf{a}\cdot\nabla\textbf{u}\right)\cdot\textbf{v}^h\,d\Omega}+\int_{\Omega}{\tilde{q}\left(\nabla\cdot\textbf{u}\right)\,d\Omega}\text{ .}
\end{equation}
Integration by parts on both terms gives
\begin{equation}
c\left(\textbf{a},\textbf{u},\textbf{v}^h\right) = -\int_{\Omega}{\textbf{u}\cdot\left(\textbf{a}\cdot\nabla\textbf{v}^h\right)d\Omega}-\int_{\Omega}{\nabla\tilde{q}\cdot\textbf{u}\,d\Omega}\text{ .}
\end{equation}
\noindent
Grouping the quantities dotted with $\mathbf{u}$ and introducing $\tau_\text{M}^{-\frac{1}{2}}\tau_\text{M}^{\frac{1}{2}} = 1$ yields
\begin{align}
 c\left(\textbf{a},\textbf{u},\textbf{v}^h\right) &= -\int_{\Omega}{\textbf{u}\cdot\left(\textbf{a}\cdot\nabla\textbf{v}^h+\nabla\tilde{q}\right)d\Omega} \\
&=-\int_{\Omega}{\tau_\text{M}^{-\frac{1}{2}}\textbf{u}\cdot\tau_\text{M}^{\frac{1}{2}}\left(\textbf{a}\cdot\nabla\textbf{v}^h+\nabla\tilde{q}\right)d\Omega} \\
&\le{\vertii{\tau_\text{M}^{-\frac{1}{2}}\textbf{u}}}\,{\vertii{\tau_\text{M}^{\frac{1}{2}}\left(\textbf{a}\cdot\nabla\textbf{v}^h+\nabla\tilde{q}\right)}} \\
&\le {\vertiii{\left(\textbf{u},p,\tilde{p}\right)}}_\ast{\vertiii{\left(\textbf{v}^h,q^h,\tilde{q}\right)}}\text{ .}
\end{align}
\noindent

\noindent
$\textbf{Continuity of Term 2}$: The second term of $A_\text{red}\left(\left(\textbf{u},p,\tilde{p}\right),\left(\textbf{v}^h,q^h,\tilde{q}\right)\right)$ is
\begin{align}
 k\left(\textbf{u},\textbf{v}^h\right)&=\int_{\Omega}{2\nu\nabla^s\textbf{u}:\nabla^s\textbf{v}^hd\Omega} \\
&\le\left(\sqrt{2\nu}{\vertii{\nabla^s\textbf{u}}}\right)\left(\sqrt{2\nu}{\vertii{\nabla^s\textbf{v}^h}}\right) \\
&\le\vertiii{\left(\textbf{u},p,\tilde{p}\right)}_\ast\vertiii{\left(\textbf{v}^h,q^h,\tilde{q}\right)}\text{ .}
\end{align}

\noindent
$\textbf{Continuity of Term 3}$: The third term is
\begin{align}
-b\left(\textbf{v}^h,p\right) &=-\int_{\Omega}{\nabla\cdot\textbf{v}^h p\,d\Omega}\\
&\le \left|\int_{\Omega}\left(\tau_\text{C}^{-\frac{1}{2}}p\right)\left(\tau_\text{C}^{\frac{1}{2}}\nabla\cdot\textbf{v}^h\right)\,d\Omega\right| \\
 &\le {\vertii{\tau_\text{C}^{-\frac{1}{2}}p}}{\vertii{\tau_\text{C}^{\frac{1}{2}}\left(\nabla\cdot\textbf{v}^h\right)}} \\
 &\le \vertiii{\left(\textbf{u},p,\tilde{p}\right)}_\ast\vertiii{\left(\textbf{v}^h,q^h,\tilde{q}\right)}\text{ .}
\end{align}

\noindent
$\textbf{Continuity of Term 4}$: The fourth term, $b\left(\mathbf{u},q^h\right)$, is zero because $\mathbf{u}\in\mathring{\mathcal{V}}_{\mathcal{Q}^h}$.

\noindent
$\textbf{Continuity of Term 5}$: The fifth term can be split into two parts, 
\begin{align}
\nonumber \left(\textbf{a}\cdot\nabla\textbf{u}-\nabla\cdot\left(2\nu\nabla^s\textbf{u}\right)+\nabla\tilde{p},\tau_\text{M}\left(\textbf{a}\cdot\nabla\textbf{v}^h+\nabla\tilde{q}\right)\right) =~& \left(\textbf{a}\cdot\nabla\textbf{u}+\nabla\tilde{p},\tau_\text{M}\left(\textbf{a}\cdot\nabla\textbf{v}^h+\nabla\tilde{q}\right)\right) \\
&-\left(\nabla\cdot\left(2\nu\nabla^s\textbf{u}\right),\tau_\text{M}\left(\textbf{a}\cdot\nabla\textbf{v}^h+\nabla\tilde{q}\right)\right)\text{ .}\label{eq:split-term-5}
\end{align}
\noindent
The first term on the right-hand side of \eqref{eq:split-term-5} is bounded as follows:
\begin{align}
\left(\textbf{a}\cdot\nabla\textbf{u}+\nabla\tilde{p},\tau_\text{M}\left(\textbf{a}\cdot\nabla\textbf{v}^h+\nabla\tilde{q}\right)\right) &= \int_{\Omega}{\left(\textbf{a}\cdot\nabla\textbf{u}+\nabla\tilde{p}\right)\cdot\left(\tau_\text{M}\left(\textbf{a}\cdot\nabla\textbf{v}^h+\nabla\tilde{q}\right)\right)\,d\Omega} \\
 &\le {\vertii{\tau_\text{M}^{\frac{1}{2}}\left(\textbf{a}\cdot\nabla\textbf{u}+\nabla\tilde{p}\right)}}\,{\vertii{\tau_\text{M}^{\frac{1}{2}}\left(\textbf{a}\cdot\nabla\textbf{v}^h+\nabla\tilde{q}\right)}} \\
 &\le \vertiii{\left(\textbf{u},p,\tilde{p}\right)}_\ast\vertiii{\left(\textbf{v}^h,q^h,\tilde{q}\right)}\text{ .}
\end{align}
\noindent
The second term of \eqref{eq:split-term-5} is bounded above by:
\begin{align}
-\left(\nabla\cdot\left(2\nu\nabla^s\textbf{u}\right),\tau_\text{M}\left(\textbf{a}\cdot\nabla\textbf{v}^h+\nabla\tilde{q}\right)\right) & \le{\vertii{\tau_\text{M}^{\frac{1}{2}}\left(\nabla\cdot\left(2\nu\nabla^s\textbf{u}\right)\right)}}\,{\vertii{\tau_\text{M}^{\frac{1}{2}}\left(\textbf{a}\cdot\nabla\textbf{v}^h+\nabla\tilde{q}\right)}} \\
 & \le\vertiii{\left(\textbf{u},p,\tilde{p}\right)}_\ast\vertiii{\left(\textbf{v}^h,q^h,\tilde{q}\right)}\text{ .}
\end{align}
\noindent
Thus,
\begin{equation}
\left(\textbf{a}\cdot\nabla\textbf{u}-\nabla\cdot\left(2\nu\nabla^s\textbf{u}\right)+\nabla\tilde{p},\tau_\text{M}\left(\textbf{a}\cdot\nabla\textbf{v}^h+\nabla\tilde{q}\right)\right)\le2\vertiii{\left(\textbf{u},p,\tilde{p}\right)}_\ast\vertiii{\left(\textbf{v}^h,q^h,\tilde{q}\right)}\text{ .}
\end{equation}
\noindent

\noindent
$\textbf{Continuity of Term 6}$: Lastly, the term
\begin{equation}
\left(\tau_\text{C}\nabla\cdot\textbf{u},\nabla\cdot\textbf{v}^h\right)= \int_{\Omega}\left(\tau_\text{C}\nabla\cdot\textbf{u}\right)\left(\nabla\cdot\textbf{v}^h\right)d\Omega
\end{equation}
can be bounded above by:
\begin{align}
\left(\tau_\text{C}\nabla\cdot\textbf{u},\nabla\cdot\textbf{v}^h\right) &\le {\vertii{\tau_\text{C}^{\frac{1}{2}}\nabla\cdot\textbf{u}}}\,{\vertii{\tau_\text{C}^{\frac{1}{2}}\nabla\cdot\textbf{v}^h}} \\
 &\le \vertiii{\left(\textbf{u},p,\tilde{p}\right)}_\ast\vertiii{\left(\textbf{v}^h,q^h,\tilde{q}\right)}\text{ .}
\end{align}

\noindent\newline
Combining the bounds on terms 1--6, we obtain \eqref{eq:A-bdd}.
\end{proof}
\end{lemma}

We continue our convergence analysis by proving inf-sup stability of the bilinear form $A_\text{red}$.
\begin{lemma}[Inf-Sup Stability]\label{Inf-Sup Stability}
There exists a constant $\gamma > 0$ independent of $h,\nu, \textbf{a}$ such that:
\begin{equation}
\inf_{(\mathbf{u}^h,p^h,\tilde{p})\in\mathcal{V}^h\times\mathcal{Q}^h\times\tilde{\mathcal{Q}}}~\sup_{(\mathbf{v}^h,q^h,\tilde{q})\in\mathcal{V}^h\times\mathcal{Q}^h\times\tilde{\mathcal{Q}}}\frac{A_\text{red}\left((\mathbf{u}^h,p^h,\tilde{p}),(\mathbf{v}^h,q^h,\tilde{q})\right)}{\vertiii{(\mathbf{u}^h,p^h,\tilde{p})}\vertiii{(\mathbf{v}^h,q^h,\tilde{q})}}\geq \gamma\text{ .} \label{Inf-Sup-Final}
\end{equation}

\begin{proof}

% [SC] 

Let $\left(\textbf{u}^h,p^h,\tilde{p}\right)\in\mathcal{V}^h\times\mathcal{Q}^h\times\tilde{\mathcal{Q}}$.\footnote{N.b. that this is an arbitrary function in the stated spaces, and not necessarily the solution to the discrete problem.} Because $\nabla\cdot\mathbf{a}=0$, $c(\mathbf{a},\mathbf{u}^h,\mathbf{u}^h) = c_\text{skew}(\mathbf{a},\mathbf{u}^h,\mathbf{u}^h) = 0$, so
\begin{equation}
\begin{split}
A_{\text{red}}\left(\left(\textbf{u}^h,p^h,\tilde{p}\right),\left(\textbf{u}^h,p^h,\tilde{p}\right)\right) =~& k\left(\textbf{u}^h,\textbf{u}^h\right)+\left(\tau_\text{M}\left(\textbf{a}\cdot\nabla\textbf{u}^h+\nabla\tilde{p}\right),\textbf{a}\cdot\nabla\textbf{u}^h+\nabla\tilde{p}\right) \\
&+ \left(\tau_\text{M}\left(\textbf{a}\cdot\nabla\textbf{u}^h+\nabla\tilde{p}\right),\nabla\cdot\left(2\nu\nabla^s\textbf{u}^h\right)\right) \\
&+\left(\tau_\text{C}\nabla\cdot\textbf{u}^h,\nabla\cdot\textbf{u}^h\right)\text{ .}
\end{split}
\end{equation}
From \eqref{eq:stab-param-defn} and \eqref{eq:inv-est}, it holds that
\begin{equation}
\begin{split}
A_{\text{red}}\left(\left(\textbf{u}^h,p^h,\tilde{p}\right),\left(\textbf{u}^h,p^h,\tilde{p}\right)\right) \geq~& \frac{1}{2}k\left(\textbf{u}^h,\textbf{u}^h\right)+\frac{1}{2}\left(\tau_\text{M}\left(\textbf{a}\cdot\nabla\textbf{u}^h+\nabla\tilde{p}\right),\textbf{a}\cdot\nabla\textbf{u}^h+\nabla\tilde{p}\right) \\
&+ \left(\tau_\text{C}\nabla\cdot\textbf{u}^h,\nabla\cdot\textbf{u}^h\right)
\end{split}
\end{equation}
and
\begin{equation}
k\left(\textbf{u}^h,\textbf{u}^h\right)\geq\left(\tau_\text{M}\nabla\cdot\left(2\nu\nabla^s\textbf{u}^h\right),\ \nabla\cdot\left(2\nu\nabla^s\textbf{u}^h\right)\right)\text{ .}
\end{equation}
Thus:
\begin{equation}
\begin{split}
A_{\text{red}}\left(\left(\textbf{u}^h,p^h,\tilde{p}\right),\left(\textbf{u}^h,p^h,\tilde{p}\right)\right) \geq~&
\frac{1}{4}k\left(\textbf{u}^h,\textbf{u}^h\right)+\frac{1}{2}\left(\tau_\text{M}\left(\textbf{a}\cdot\nabla\textbf{u}^h+\nabla\tilde{p}\right),\textbf{a}\cdot\nabla\textbf{u}^h+\nabla\tilde{p}\right) \\
&+ \frac{1}{4}\left(\tau_\text{M}\nabla\cdot\left(2\nu\nabla^s\textbf{u}^h\right),\ \nabla\cdot\left(2\nu\nabla^s\textbf{u}^h\right)\right) \\
&+ \left(\tau_\text{C}\nabla\cdot\textbf{u}^h,\nabla\cdot\textbf{u}^h\right)
\end{split}
\end{equation}
Applying Assumption \ref{assumption1} with $q^h = p^h$ implies that there exists $\mathbf{w}^h\in\mathcal{V}^h$ such that
\begin{equation}
\begin{split}
A_{\text{red}}\left(\left(\textbf{u}^h,p^h,\tilde{p}\right),\left(\textbf{w}^h,0,0\right)\right) \geq~& 
\tilde{\gamma}{\vertii{p^h}}^2+c\left(\textbf{a},\textbf{u}^h,\textbf{w}^h\right)+k\left(\textbf{u}^h,\textbf{w}^h\right) \\
&+ \left(\textbf{a}\cdot\nabla\textbf{u}^h-\nabla\cdot\left(2\nu\nabla^s\textbf{u}^h\right)+\nabla\tilde{p},\tau_\text{M}\left(\textbf{a}\cdot\nabla\textbf{w}^h\right)\right) \\
&+\left(\tau_\text{C}\nabla\cdot\textbf{u}^h,\nabla\cdot\textbf{w}^h\right)\text{ .} \label{A_red_wh_0_0}
\end{split}
\end{equation}
Using Young's and Poincar\'{e}'s inequalities, we obtain the following bounds on the convective and diffusive terms:
\begin{equation}
c\left(\textbf{a},\textbf{u}^h,\textbf{w}^h\right)\le\frac{5}{2\tilde{\gamma}}\left(C_\text{Poin}^2L^2{\vertii{\textbf{a}}}_{L^\infty}^2\nu^{-1}\right)k\left(\textbf{u}^h,\textbf{u}^h\right)+\frac{\tilde{\gamma}}{10}{\vertii{\nabla\textbf{w}^h}}^2
\end{equation}
\begin{equation}
k\left(\textbf{u}^h,\textbf{w}^h\right)\le\frac{5}{2\tilde{\gamma}}\nu k\left(\textbf{u}^h,\textbf{u}^h\right)+\frac{\tilde{\gamma}}{10}{\vertii{\nabla\textbf{w}^h}}^2
\end{equation}
where $C_\text{Poin}>0$ is the $h$-independent dimensionless constant from Poincar\'{e}'s inequality and $L$ is a length scale proportional to the diameter of $\Omega$. Using Young's inequality,
\begin{equation}
\left(\tau_\text{C}\nabla\cdot\textbf{u}^h,\nabla\cdot\textbf{w}^h\right)\le\frac{5}{2\tilde{\gamma}}\tau_\text{C}\left(\tau_\text{C}\nabla\cdot\textbf{u}^h,\nabla\cdot\textbf{u}^h\right)+\frac{\tilde{\gamma}}{10}{\vertii{\nabla\textbf{w}^h}}^2\text{ .}
\end{equation}
Finally with \eqref{eq:stab-param-defn}, we obtain the following bounds:
\begin{equation}
\left(\textbf{a}\cdot\nabla\textbf{u}^h+\nabla\tilde{p},\tau_\text{M}\left(\textbf{a}\cdot\nabla\textbf{w}^h\right)\right)\le\frac{5}{2\tilde{\gamma}}\left(\frac{h{\vertii{\textbf{a}}}_{L^\infty}}{2}\right)\left(\tau_\text{M}\left(\textbf{a}\cdot\nabla\textbf{u}^h+\nabla\tilde{p}\right),\textbf{a}\cdot\nabla\textbf{u}^h+\nabla\tilde{p}\right)+\frac{\tilde{\gamma}}{10}{\vertii{\nabla\textbf{w}^h}}^2
\end{equation}
and
\begin{equation}
\left(\nabla\cdot\left(2\nu\nabla^s\textbf{u}^h\right),\tau_\text{M}\left(\textbf{a}\cdot\nabla\textbf{w}^h\right)\right)\le\frac{5}{2\tilde{\gamma}}\left(\frac{h{\vertii{\textbf{a}}}_{L^\infty}}{2}\right)\left(\tau_\text{M}\nabla\cdot\left(2\nu\nabla^s\textbf{u}^h\right),\nabla\cdot\left(2\nu\nabla^s\textbf{u}^h\right)\right)+\frac{\tilde{\gamma}}{10}{\vertii{\nabla\textbf{w}^h}}^2\text{ .}
\end{equation}
Combining these results, \eqref{A_red_wh_0_0} is bounded below by:
\begin{equation}
\begin{split}
A_{\text{red}}\left(\left(\textbf{u}^h,p^h,\tilde{p}\right),\left(\textbf{w}^h,0,0\right)\right) \geq~& 
\frac{\tilde{\gamma}}{10}{\vertii{p^h}}^2-\frac{5}{2\tilde{\gamma}}\left(C_{Poin}^2L^2{\vertii{\textbf{a}}}_{L^\infty}^2\nu^{-1}\right)k\left(\textbf{u}^h,\textbf{u}^h\right)-\frac{5}{2\tilde{\gamma}}\nu k\left(\textbf{u}^h,\textbf{u}^h\right) \\
&- \frac{5}{2\tilde{\gamma}}\left(\frac{h{\vertii{\textbf{a}}}_{L^\infty}}{2}\right)\left(\tau_\text{M}\left(\textbf{a}\cdot\nabla\textbf{u}^h+\nabla\tilde{p}\right),\textbf{a}\cdot\nabla\textbf{u}^h+\nabla\tilde{p}\right) \\
&- \frac{5}{2\tilde{\gamma}}\left(\frac{h{\vertii{\textbf{a}}}_{L^\infty}}{2}\right)\left(\tau_\text{M}\nabla\cdot\left(2\nu\nabla^s\textbf{u}^h\right),\nabla\cdot\left(2\nu\nabla^s\textbf{u}^h\right)\right) \\
&- \frac{5}{2\tilde{\gamma}}\tau_\text{C}\left(\tau_\text{C}\nabla\cdot\textbf{u}^h,\nabla\cdot\textbf{u}^h\right)\text{ .}
\end{split}
\end{equation}
Now let
\begin{equation}
\tilde{C} = \frac{\tilde{\gamma}}{40\nu}\min\left\{1,\text{Pe}^{-2},\text{Pe}_h^{-1}\right\}\text{ .}
\end{equation}
Then it holds that
\begin{equation}
\begin{split}
A_{\text{red}}\left(\left(\textbf{u}^h,p^h,\tilde{p}\right),\left(\textbf{u}^h+\tilde{C}\textbf{w}^h,p^h,\tilde{p}\right)\right) \geq~& \frac{1}{8}k\left(\textbf{u}^h,\textbf{u}^h\right)+\frac{1}{2}\left(\tau_\text{M}\left(\textbf{a}\cdot\nabla\textbf{u}^h+\tilde{p}\right),\textbf{a}\cdot\nabla\textbf{u}^h+\tilde{p}\right) \\
&+ \frac{1}{4}\left(\tau_\text{M}\nabla\cdot\left(2\nu\nabla^s\textbf{u}^h\right),\nabla\cdot\left(2\nu\nabla^s\textbf{u}^h\right)\right) \\
&+\frac{1}{2}\left(\tau_\text{C}\nabla\cdot\textbf{u}^h,\nabla\cdot\textbf{u}^h\right) \\
&+\frac{\tilde{\gamma}^2}{400}\nu^{-1}\min{\left\{1,{\rm Pe}^{-2},{\rm Pe}_h^{-1}\right\}}{\vertii{p^h}}^2\text{ .}
\end{split}
\end{equation}
Recalling \eqref{triple-bar-norm} and \eqref{Define:alpha}, we have
\begin{equation}
A_{\text{red}}\left(\left(\textbf{u}^h,p^h,\tilde{p}\right),\left(\textbf{u}^h+\tilde{C}\textbf{w}^h,p^h,\tilde{p}\right)\right)\geq\beta{\vertiii{\left(\textbf{u}^h,p^h,\tilde{p}\right)}}^2\text{ ,}
\end{equation}
where
\begin{equation}
\beta=\min{\left\{\frac{1}{8},\frac{{\tilde{\gamma}}^2}{400}\right\}}\text{ .}
\end{equation}
From the triangle inequality, we have:
\begin{equation}
\vertiii{\left(\textbf{u}^h+\tilde{C}\textbf{w}^h,p^h,\tilde{p}\right)} \le \vertiii{\left(\textbf{u}^h,p^h,\tilde{p}\right)}+\tilde{C}\vertiii{\left(\textbf{w}^h,0,0\right)}\text{ .}
\end{equation}
From here, to obtain the desired inf-sup stability result given by \eqref{Inf-Sup-Final}, we must first prove the existence of a constant $\xi$ independent of $h,\nu,$ and $\textbf{a}$ such that:
\begin{equation}
{\tilde{C}}^2{\vertiii{\left(\textbf{w}^h,0,0\right)}}^2\le\xi\alpha{\vertii{p^h}}^2\text{ .} \label{eq.existence-xi-result}
\end{equation}
The constant $\xi$ can be obtained by substituting $\left({\textbf{w}}^h,0,0\right)$ into the norm defined by \eqref{triple-bar-norm} and simplifying the remaining terms:
\begin{align}
\nonumber \vertiii{({\textbf{w}}^h,0,0)}^2 =~& \left(2\nu\nabla^s{\textbf{w}}^h,\nabla^s{\textbf{w}}^h\right)+\left(\tau_\text{C}\nabla\cdot{\textbf{w}}^h,\nabla\cdot{\textbf{w}}^h\right) \\
&+\left(\tau_\text{M}\left(\textbf{a}\cdot\nabla{\textbf{w}}^h\right),\textbf{a}\cdot\nabla{\textbf{w}}^h\right)+\left(\tau_\text{M}\nabla\cdot\left(2\nu\nabla^s{\textbf{w}}^h\right),\nabla\cdot\left(2\nu\nabla^s{\textbf{w}}^h\right)\right)\text{ .} \label{eq.lineabove}
\end{align}
Using Assumption \ref{assumption2},
\begin{align}
\nonumber\vertiii{\left(\textbf{w}^h,0,0\right)}^2\leq&\left(2\nu\nabla^s\textbf{w}^h,\nabla^s\textbf{w}^h\right)+\left(\tau_\text{C}\nabla\cdot\textbf{w}^h,\nabla\cdot\textbf{w}^h\right)\\
&+\left(\tau_\text{M}\left(\textbf{a}\cdot\nabla\textbf{w}^h\right),\textbf{a}\cdot\nabla\textbf{w}^h\right)+\tau_\text{M}\frac{C_\text{inv}\nu^2}{h^2}\left(\nabla^s\textbf{w}^h,\nabla^s\textbf{w}^h\right)\text{ .} \label{eq.sub-younginequality-result-here}
\end{align}
Recalling the definition of $\mathbf{w}^h$, we have
\begin{align}
    \tilde{C}^2\vertiii{(\mathbf{w}^h,0,0)}^2 \leq \left(\frac{\tilde{\gamma}}{40}\right)^2\alpha^2\left(2\nu + \tau_\text{C} + \tau_\text{M}\vert\mathbf{a}\vert^2 + \tau_\text{M}\frac{C_\text{inv}\nu^2}{h^2}\right)\left\Vert p^h\right\Vert^2\text{ .}
\end{align}
By choosing convenient branches from the minimums in $\alpha$ and $\tau_\text{M}$ to obtain upper bounds, and considering both branches of the maximum in $\tau_\text{C}$, we easily obtain
\begin{equation}
    2\nu\alpha \leq 2~,\quad \alpha\tau_\text{C}\leq 2~,\quad\alpha\tau_\text{M}\vert\mathbf{a}\vert^2\leq 1~,\quad\text{and}\quad \frac{\alpha \tau_\text{M}C_\text{inv}\nu^2}{h^2}\leq 1\text{ ,}
\end{equation}
which establishes \eqref{eq.existence-xi-result}, with
\begin{equation}
    \xi = 6\left(\frac{\tilde{\gamma}}{40}\right)^2\text{ .}
\end{equation} It then holds that there exists a constant $\gamma$, independent of $h,\nu,$ and $\textbf{a}$, such that:
\begin{equation}
    \inf_{(\mathbf{u}^h,p^h,\tilde{p})\in\mathcal{V}^h\times\mathcal{Q}^h\times\tilde{\mathcal{Q}}}~\sup_{(\mathbf{v}^h,q^h,\tilde{q})\in\mathcal{V}^h\times\mathcal{Q}^h\times\tilde{\mathcal{Q}}}\frac{A_\text{red}\left((\mathbf{u}^h,p^h,\tilde{p}),(\mathbf{v}^h,q^h,\tilde{q})\right)}{\vertiii{(\mathbf{u}^h,p^h,\tilde{p})}\vertiii{(\mathbf{v}^h,q^h,\tilde{q})}}\geq \gamma\text{ ,}
\end{equation}
where $\gamma$ has been found to be:
\begin{equation}
\gamma = \frac{\beta}{\left(1+\xi^\frac{1}{2}\right)}\text{ .}
\end{equation}
\end{proof}
\end{lemma}

The preceding consistency (Lemma \ref{Consistency}), continuity (Lemma \ref{Continuity}), and stability (Lemma \ref{Inf-Sup Stability}) results combine in a standard way to bound the discrete solution's error.
\begin{theorem}[Error bound]\label{thm:err-bnd}
Suppose that the exact solution $(\mathbf{u},p,p)\in (H^2(\Omega))^d\times H^1(\Omega)\times H^1(\Omega)$.  Then 
\begin{equation}
    \vertiii{\left(\mathbf{u}-\mathbf{u}^h,p-p^h,p-\tilde{p}\right)} \leq C\inf_{\substack{\mathbf{w}^h\in\mathring{\mathcal{V}}_{\mathcal{Q}^h}\cap\mathcal{V}^h\\r^h,s^h\in\mathcal{Q}^h}}\vertiii{\left(\mathbf{u}-\mathbf{w}^h,p-r^h,p-s^h\right)}_*\text{ ,}\label{eq:error-bound}
\end{equation}
where $C > 0$ is a constant independent of $h$, $\nu$, and $\mathbf{a}$.
%
%[DK] Probably "standard" enough to simply state without proof.  (Maybe cite a textbook or something.)
%
%[DK] The restriction to discretely div-free $\mathbf{w}^h$ in the infimum follows from the corresponding restriction on the trial function slot in the continuity lemma, which was needed to add a convenient term in the proof (under "Continuity of Term 1").
%
\end{theorem}

To determine the convergence of the right side of \eqref{eq:error-bound} to zero with respect to discretization refinement, we must introduce an assumption regarding the approximation power of $\mathcal{V}^h$ and $\mathcal{Q}^h$.  In particular, we consider the approximation power of discretely divergence-free functions in $\mathcal{V}^h$, i.e., the subspace $\mathring{\mathcal{V}}_{\mathcal{Q}^h}\cap\mathcal{V}^h$ over which the infimum is taken in \eqref{eq:error-bound}.  The best interpolation error in this subspace can be bounded by the error in the corresponding Stokes projector.  (Assuming $\mathcal{V}^h$ and $\mathcal{Q}^h$ are an optimal and inf-sup stable velocity--pressure pair for the Stokes problem, the Stokes projector's error should converge optimally.)
\begin{assumption}\label{assmpt:approx}
There exist $k_v,k_q\in\mathbb{Z}_{\geq 0}$, $C_\text{int}\in\mathbb{R}_{> 0}$ independent of the global quasiuniform element size $h$, and interpolation operators $\mathcal{I}_v : \mathring{\mathcal{V}} \to \mathring{\mathcal{V}}_{\mathcal{Q}^h}\cap\mathcal{V}^h$ and $\mathcal{I}_q : \mathcal{Q} \to \mathcal{Q}^h$, such that, $\forall\mathbf{v} \in \mathring{\mathcal{V}} \cap (H^{k_v+1}(\Omega))^d$ and nonnegative integer $l \leq k_v + 1$,
\begin{equation}
\sum_{e=1}^{n_\text{el}} | \mathbf{v} - \mathcal{I}_v \mathbf{v} |^2_{(H^l(\Omega^e))^d} \leq C_\text{int} h^{2(k_v+1-l)} | \mathbf{v} |^2_{(H^{k_v+1}(\Omega))^d}\text{ ,}\label{eq:v-approx}
\end{equation}
and, $\forall q \in \mathcal{Q} \cap H^{k_q+1}(\Omega)$ and nonnegative integer $l \leq k_q + 1$,
\begin{equation}
\sum_{e=1}^{n_\text{el}} | q - \mathcal{I}_q q |^2_{H^l(\Omega^e)} \leq C_\text{int} h^{2(k_q+1-l)} | q |^2_{H^{k_q+1}(\Omega)}\text{ ,}\label{eq:q-approx}
\end{equation}
where $\mathring{\mathcal{V}}$ is the solenoidal subspace of $\mathcal{V}$.  The integers $k_v$ and $k_q$ are interpreted as the degrees of the discrete velocity and pressure spaces.
\end{assumption}

\begin{theorem}[Convergence]\label{thm:conv}
Given the preceding assumptions, we have that, for an exact solution $(\mathbf{u},p,p)\in (H^{k_v+1}(\Omega))^d\times H^{k_q+1}(\Omega)\times H^{k_q+1}(\Omega)$,
\begin{equation}
    \vertiii{\left(\mathbf{u}-\mathbf{u}^h,p-p^h,p-\tilde{p}\right)}^2 \leq C\left(\max\{h\vert\mathbf{a}\vert,\nu\}h^{2k_v}\vert\mathbf{u}\vert^2_{(H^{k_v+1}(\Omega))^d} + \frac{1}{\max\{h\vert\mathbf{a}\vert,\nu\}}h^{2k_q+2}\vert p\vert^2_{(H^{k_q+1}(\Omega))^d}\right)\text{ ,}\label{eq:conv-result}
\end{equation}
where $C\in\mathbb{R}_{> 0}$ is independent of $h$, $\nu$, and $\mathbf{a}$.
\begin{proof}
Invoking Theorem \ref{thm:err-bnd}, we can bound the total error by the interpolation error,
\begin{align}
\nonumber\vertiii{\left(\bm{\eta}_u,\eta_p,\eta_p\right)}_*^2 =~& k\left(\bm{\eta}_u,\bm{\eta}_u\right)+\alpha\left\Vert \eta_p\right\Vert^2 + \left\Vert\tau_\text{C}^{1/2}\nabla\cdot\bm{\eta}_u\right\Vert^2\\
\nonumber &+ \left\Vert\tau_\text{M}^{1/2}\left(\mathbf{a}\cdot\nabla\bm{\eta}_u + \nabla\eta_p\right)\right\Vert^2 + \left\Vert\tau_\text{M}^{1/2}\nabla\cdot\left(2\nu\nabla^s\bm{\eta}_u\right)\right\Vert^2 \\
&+ \left\Vert\tau_\text{M}^{-1/2}\bm{\eta}_u\right\Vert^2 + \left\Vert\tau_\text{C}^{-1/2}\eta_p\right\Vert^2\text{ ,}\label{eq:full-interp-err}
\end{align}
where
\begin{equation}
    \bm{\eta}_u := \mathbf{u} - \mathcal{I}_v\mathbf{u}\quad\text{and}\quad\eta_p := p - \mathcal{I}_p p\text{ .}
\end{equation}
We now bound each term of the right-hand side of \eqref{eq:full-interp-err}.
In the remainder of this proof, we use ``$C$'' to represent a generic constant independent of $h$, $\vert\mathbf{a}\vert$, and $\nu$; its numerical value may be different in different places.
A bound on the first term follows immediately from \eqref{eq:v-approx} and Korn's inequality:
\begin{equation}
    k(\bm{\eta}_u,\bm{\eta}_u) \leq C\nu h^{2k_v}\vert\mathbf{u}\vert^2_{(H^{k_v+1}(\Omega))^d}\text{ .}
\end{equation}
Using \eqref{Define:alpha} and \eqref{eq:q-approx}, the second term is bounded as
\begin{equation}
    \alpha\Vert\eta_p\Vert^2_{L^2(\Omega)} \leq \frac{C}{\max\{h\vert\mathbf{a},\nu\}} h^{2k_q+2}\vert p\vert_{H^{k_q+1}(\Omega)}^2\text{ .}
\end{equation}
A bound on the third term follows from \eqref{eq:stab-param-defn}, \eqref{eq:v-approx}, and bounding the divergence by the gradient:
\begin{equation}
    \left\Vert\tau_\text{C}^{1/2}\nabla\cdot\bm{\eta}_u\right\Vert_{L^2(\Omega)}^2 \leq C\max\{h\vert\mathbf{a}\vert,\nu\} h^{2k_v}\vert\mathbf{u}\vert^2_{(H^{k_v+1}(\Omega))^d}\text{ .}
\end{equation}
We can split the bound for the fourth term on the right of \eqref{eq:full-interp-err} into two parts,
\begin{equation}\label{eq:split-bound-4th-term}
    \left\Vert\tau_\text{M}^{1/2}\left(\mathbf{a}\cdot\nabla\bm{\eta}_u + \nabla\eta_p\right)\right\Vert^2 \leq 2\left(\left\Vert \tau_\text{M}^{1/2}\mathbf{a}\cdot\nabla\bm{\eta}_u\right\Vert^2 + \left\Vert\tau_\text{M}^{1/2}\nabla\eta_p\right\Vert^2\right)\text{ ,} 
\end{equation}
where Young's inequality bounds the $(\tau_\text{M}\mathbf{a}\cdot\nabla\bm{\eta}_u,\nabla\eta_p)$ cross terms.  The terms on the right of \eqref{eq:split-bound-4th-term} can in turn be bounded using interpolation estimates of Assumption \ref{assmpt:approx} and the definition of $\tau_\text{M}$ in \eqref{eq:stab-param-defn}:
\begin{equation}
    \left\Vert \tau_\text{M}^{1/2}\mathbf{a}\cdot\nabla\bm{\eta}_u\right\Vert^2 \leq C\vert\mathbf{a}\vert h^{2k_v+1}\vert\mathbf{u}\vert^2_{(H^{k_v+1}(\Omega))^d}
\end{equation}
and
\begin{equation}
    \left\Vert\tau_\text{M}^{1/2}\nabla\eta_p\right\Vert^2 \leq \frac{C}{\max\{h\vert\mathbf{a}\vert,\nu\}}h^{2k_q+2}\vert p\vert^2_{(H^{k_q+1}(\Omega))^d}\text{ .}
\end{equation}
We now invoke the inverse estimate of Assumption \ref{assumption2} and \eqref{eq:stab-param-defn} to bound the fifth term:
\begin{equation}
    \left\Vert\tau_\text{M}^{1/2}\nabla\cdot(2\nu\nabla^s\bm{\eta}_u)\right\Vert^2 \leq C\nu h^{2k_v}\vert\mathbf{u}\vert^2_{(H^{k_v+1}(\Omega))^d}\text{ .}
\end{equation}
Finally, the last two terms are bounded straightforwardly, from Assumption \ref{assmpt:approx} and \eqref{eq:stab-param-defn}:
\begin{equation}
    \left\Vert\tau_\text{M}^{-1/2}\bm{\eta}_u\right\Vert^2 \leq C\max\{h\vert\mathbf{a}\vert,\nu\}h^{2k_v}\vert\mathbf{u}\vert^2_{(H^{k_v+1}(\Omega))^d}
\end{equation}
and
\begin{equation}
    \left\Vert\tau_\text{C}^{-1/2}\eta_p\right\Vert^2 \leq \frac{C}{\max\{h\vert\mathbf{a}\vert,\nu\}}h^{2k_q+2}\vert p\vert^2_{(H^{k_q+1}(\Omega))^d}\text{ .}
\end{equation}
Summing these bounds and grouping terms results in \eqref{eq:conv-result}.
\end{proof}
\end{theorem}

The bound of Theorem \ref{thm:conv} shows the same structure as \cite[Theorem 2]{Evans2019a}, but the generalization to arbitrary inf-sup stable elements affects its interpretation.  In \cite{Evans2019a}, it was suggested to choose elements with $k_v = k_q$, to balance velocity and pressure contributions to convergence rates in the advection-dominated setting (i.e., $\max\{h\vert\mathbf{a}\vert,\nu\} = h\vert\mathbf{a}\vert$), which is more often encountered in practical problems.  This is indeed the case when selecting divergence-conforming B-spline velocity and pressure spaces, such as those used in the numerical examples of \cite{Evans2019a}.  However, many of the common inf-sup stable velocity--pressure pairs to which the present work applies have $k_v = k_q + 1$ (e.g., the Taylor--Hood element).  This balances rates in the diffusion-dominated regime (i.e., $\max\{h\vert\mathbf{a}\vert,\nu\} = \nu$, which is always approached in the ultimate limit of $h\to 0$ for fixed $\mathbf{a}$ and $\nu$), but some of the velocity space's approximation power may be considered ``wasted'' in the advection-dominated regime, where the convergence rate of $\vertiii{\left(\mathbf{u}-\mathbf{u}^h,p-p^h,p-\tilde{p}\right)}^2$ is limited to $2k_q+1$, despite the velocity contributions converging at a (pre-asymptotic) rate of $2k_q+3$. 

\begin{remark}\label{rem:different-tau_M-asymptotics}
It is possible that error contributions could be balanced for $k_v > k_q$ in the advection-dominated regime by choosing different asymptotic behavior of $\tau_M$, but this would come at the cost of weakening the $\vertiii{\cdot}$ norm.  See \cite[Remark 5.42]{John2016} for a related discussion in the context of standard SUPG, PSPG, and grad-div stabilization.
\end{remark}

\section{Numerical tests}\label{sec:numerical}
This section first verifies the result of Theorem \ref{thm:conv} for the Oseen equation (Section \ref{sec:numerical-oseen}), then looks at its extrapolation to the Navier--Stokes problem (Section \ref{sec:numerical-ns}).  As examples of inf-sup-stable elements satisfying Assumption \ref{assumption1}, our tests use the Taylor--Hood element \cite{Taylor1973,Brezzi1991} and its isogeometric B-spline counterpart \cite{Hosseini2015}.  In short, the degree-$k$ Taylor--Hood finite element method uses continuous Lagrange spaces of polynomial degree $k$ for the velocity components and of degree $k-1$ for the pressure, 
%
%[DK] TODO: Update numerical examples to use reduced continuity of velocity field.
while its B-spline analogue uses B-spline spaces of these degrees with the same continuity.

As in \cite{Evans2019a}, we implement these numerical examples using the FEniCS \cite{Logg2012} finite element automation software, and the library tIGAr \cite{Kamensky2019} extending it to isogeometric analysis (IGA) \cite{Hughes05a,CoHuBa09}.  This allows direct translation of the variational forms defining our method into the Python-based Unified Form Language (UFL) \cite{Alnaes2014}, which is compiled \cite{Kirby2006} into low-level finite element routines.  The source code for these numerical examples is available (alongside those of \cite{Evans2019a}) in a public Git repository \cite{examples-repo}. 

\subsection{The Oseen problem}\label{sec:numerical-oseen}
\begin{figure}[!ht]\centering
\includegraphics[width=0.85\textwidth]{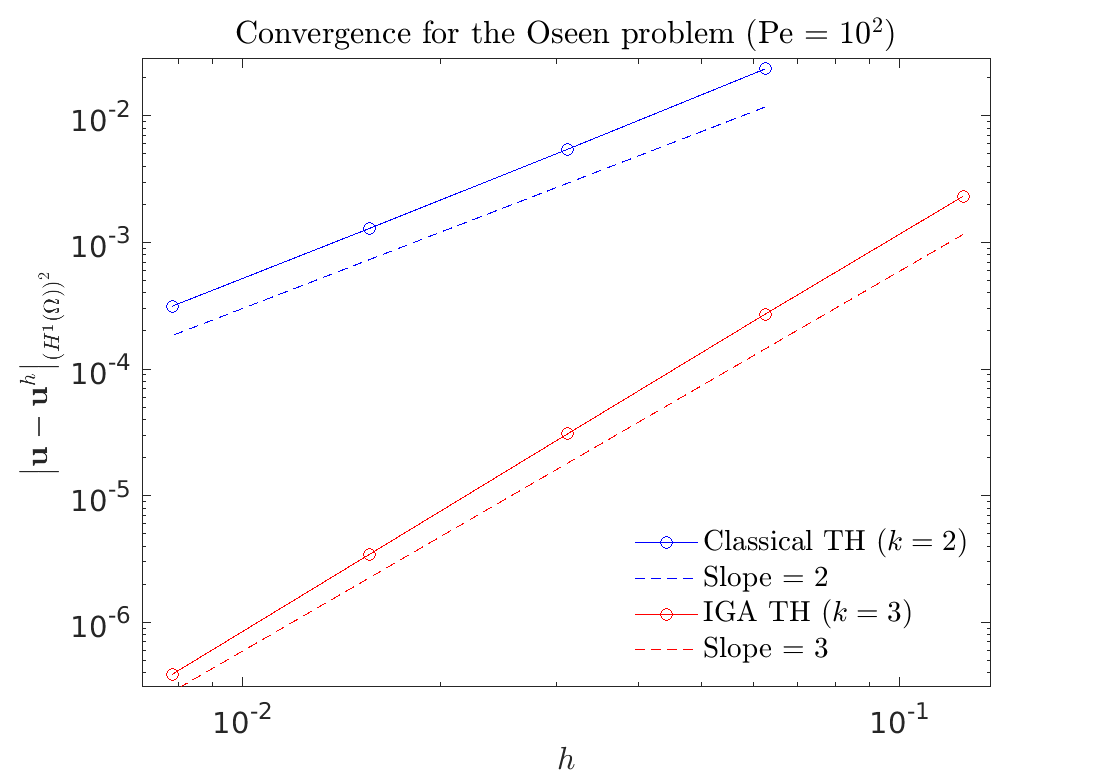}
\caption{Convergence of the $H^1$-seminorm of velocity error in the Oseen problem with $\text{Pe}=10^2$, using Taylor--Hood elements of degree $k=2$.}
\label{fig:oseen-conv-re-1e2}
\end{figure}

\begin{figure}[!ht]\centering
\includegraphics[width=0.85\textwidth]{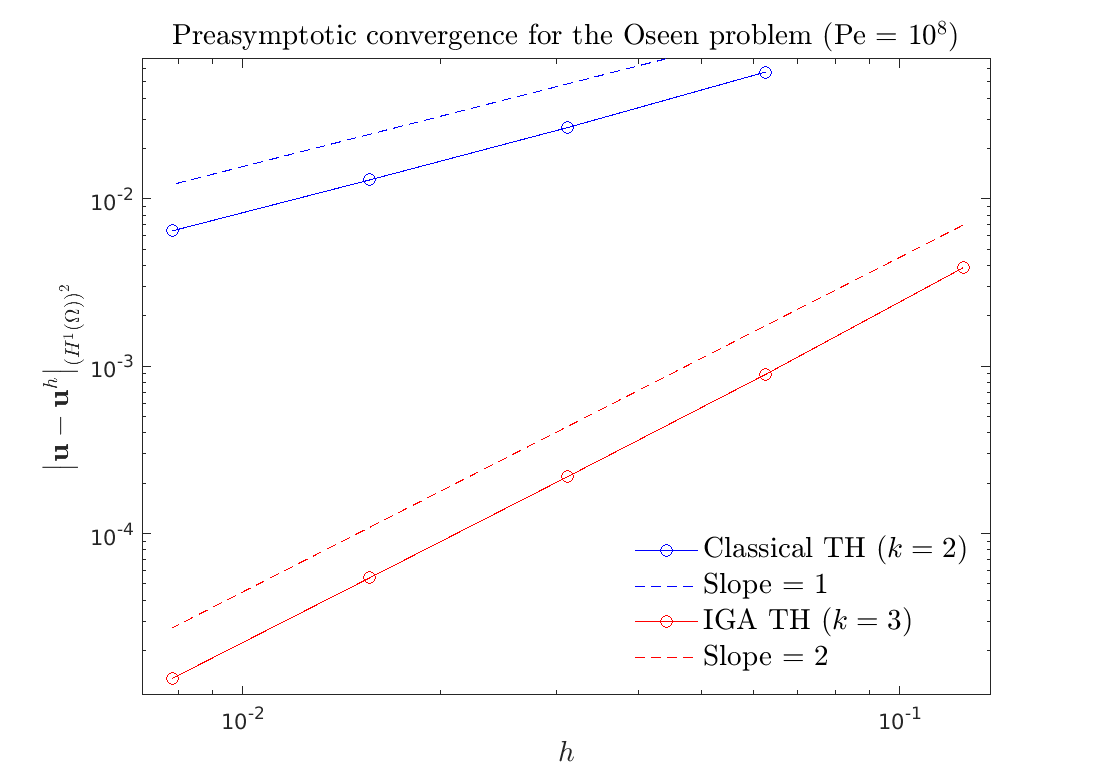}
\caption{Convergence of the $H^1$-seminorm of velocity error in the Oseen problem with $\text{Pe}=10^8$, using Taylor--Hood elements of degree $k=2$.}
\label{fig:oseen-conv-re-1e6}
\end{figure}

To verify Theorem \ref{thm:conv}, we compute solutions for the Oseen problem, using the steady flow benchmark of the regularized lid-driven cavity, which is inspired by the classic lid-driven cavity problem, but modified to have a smooth solution, so as to permit high-order convergence on quasi-uniform meshes.  Variations on the regularized lid-driven cavity have been studied previously by several groups of authors \cite{Shih1989,VanOpstal2017,Evans2019a}.  In the present work, we use the problem statement given in \cite[Section 4.2.1]{Evans2019a}, wherein the exact solution
\begin{align}
    \mathbf{u}(\mathbf{x}) &= x_1^2x_2\left(x_1^2 - 2x_1^1 + 1\right)\left(4x_2^2 - 2\right)\mathbf{e}_1 - 16x_1x_2^2\left(2x_1^2 - 3x_1 + 1\right)\left(x_2^2 - 1\right)\mathbf{e}_2\\
    p(\mathbf{x}) &= \sin(\pi x_1) \sin(\pi x_2)
\end{align}
is manufactured by an appropriate source term and Dirichlet boundary conditions on a unit square.  To match the analysis for constant $\mathbf{a}$, we set
\begin{equation}
    \mathbf{a} = \frac{\sqrt{3}}{2}\mathbf{e}_1 + \frac{1}{2}\mathbf{e}_2\text{ ,}
\end{equation}
independently of the solution $\mathbf{u}$.  We modulate the P\'{e}clet number by adjusting the kinematic viscosity $\nu$.  Because the exact solution is independent of $\nu$, this is a useful test problem to study the effect of $\text{Pe}$ on the numerical method's accuracy, while controlling for $\text{Pe}$-dependence of the exact solution (whose $H^1$ norm of velocity would often increase with $\text{Pe}$ in practical problems, due, e.g., to boundary layers becoming sharper).  Using classical Taylor--Hood elements of degree $k=2$ and isogeometric Taylor--Hood elements of degree $k=3$, we see that optimal convergence rates of $H^1$ velocity error are attained at $\text{Pe}=10^2$ (Figure \ref{fig:oseen-conv-re-1e2}).  For $\text{Pe} = 10^8$, we compare, in Figure \ref{fig:oseen-conv-re-1e6}, the preasymptotic convergence of these two elements with the rates of 1 and 2 predicted by Theorem \ref{thm:conv}, and see approximate agreement.\footnote{The theorem technically only bounds $\vertiii{(\mathbf{u}-\mathbf{u}^h,p-p^h,p-\tilde{p})}$ in a Pe-robust way, but we find that this extends to the $H^1$ norm of $\mathbf{u} - \mathbf{u}^h$ in practice.}

\begin{figure}[!ht]\centering
\includegraphics[width=0.85\textwidth]{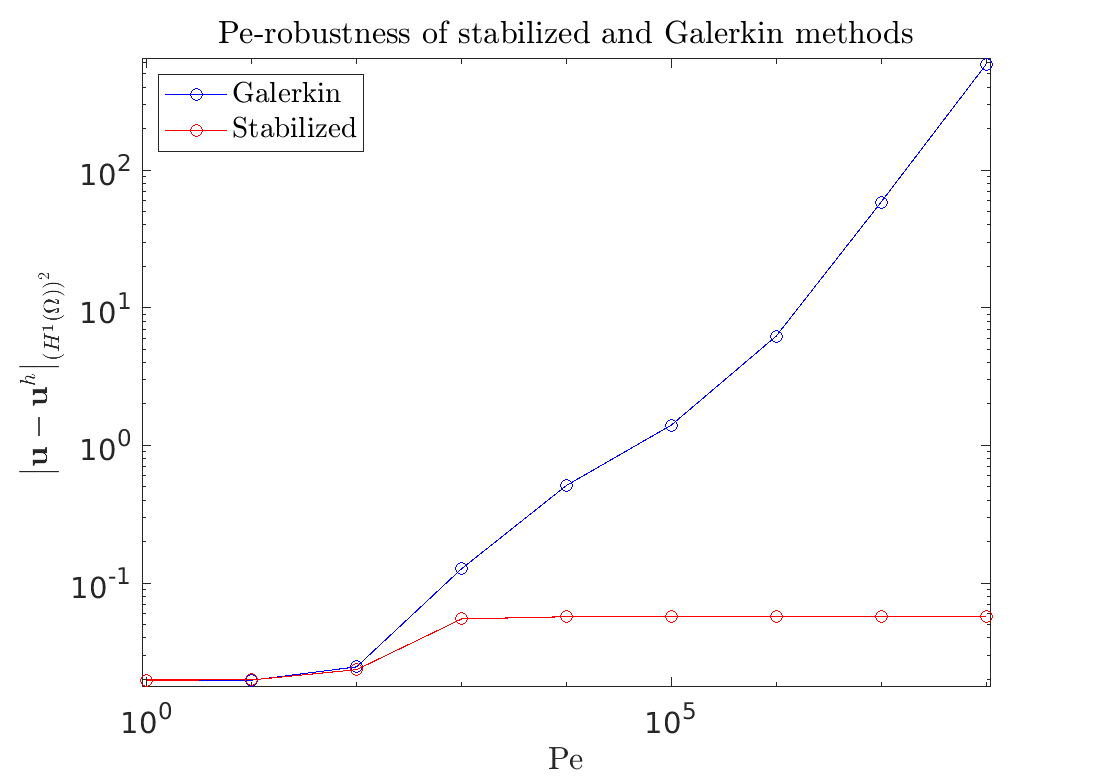}
\caption{Robustness of the proposed stabilization in the inviscid limit, using classical Taylor--Hood elements of degree $k=2$ on a fixed coarse mesh.}
\label{fig:oseen-Pe-robust}
\end{figure}

For a fixed mesh of $16\times 16$ square cells, each split into two right triangular elements, Figure \ref{fig:oseen-Pe-robust} compares the proposed method with the Galerkin method (viz., $\tau_\text{M} = \tau_\text{C} = 0$) for the classical Taylor--Hood element with $k=2$.  We see that the $H^1$ norm the stabilized scheme's velocity error reaches a plateau in the inviscid limit, while the error in the Galerkin scheme increases with $\text{Pe}$, even for a smooth, $\text{Pe}$-independent solution (which can often mask instabilities of Galerkin's method in practice).  

\subsection{The Navier--Stokes problem}\label{sec:numerical-ns}
As is often the case, the convergence result proven for Oseen flow in Section \ref{sec:oseen-conv} extends in practice to the full Navier--Stokes problem, which we demonstrate here for steady (Section \ref{sec:reg-ldc}) and unsteady (Section \ref{sec:2d-tg}) benchmarks.

\subsubsection{Steady Navier--Stokes}\label{sec:reg-ldc}
We now revisit the regularized lid-driven cavity benchmark used in Section \ref{sec:numerical-oseen}, but solve the full nonlinear Navier--Stokes problem, without fixing the advection velocity.  Numerical results for this problem were computed at Reynolds number $\text{Re}=100$.  We compare the convergence rates of the classical and isogeometric Taylor--Hood elements.  To improve nonlinear convergence, we use the same smoothed stabilization parameters from the numerical experiments of \cite{Evans2019a}:
\begin{align}
\tau_\text{M} &= \left(\textbf{u}^h\cdot\textbf{G}\cdot\textbf{u}^h+C_\text{inv}^2\nu^2\textbf{G}:\textbf{G}\right)^{-\frac{1}{2}} \label{tau_M_LDC}\\
\tau_\text{C} &= \left(\tau_\text{M} \text{tr}{\left(\textbf{G}\right)}\right)^{-1} \label{tau_C_LDC}
\end{align}
where the tensor $\textbf{G}$ is defined as
\begin{equation}
\textbf{G}={\frac{\partial \bm{\xi}}{\partial \textbf{x}}}^T\frac{\partial \bm{\xi}}{\partial \textbf{x}}\text{ ,} \label{metric-G}
\end{equation}
in which $\frac{\partial \bm{\xi}}{\partial \textbf{x}}$ is the inverse Jacobian of the mapping from the parent element to physical space.

\begin{figure}[!ht]\centering
\includegraphics[width=0.85\textwidth]{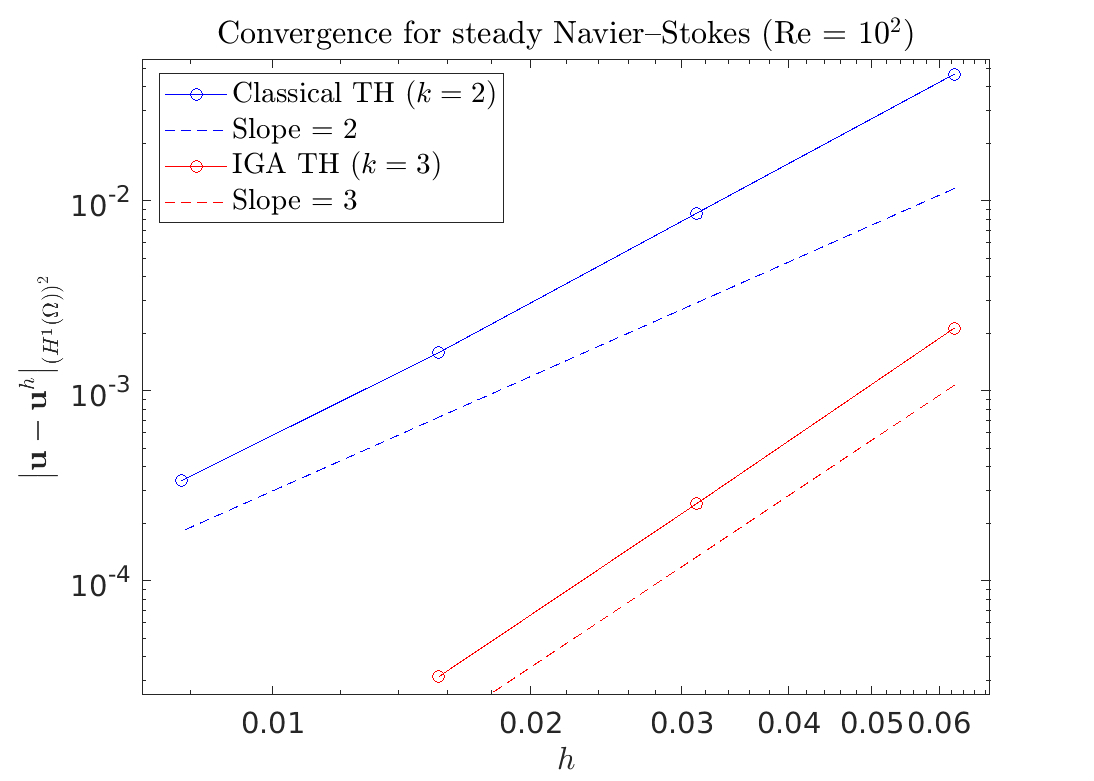}
\caption{Convergence of $H^1$ velocity error for the regularized lid-driven cavity, using classical and isogeometric Taylor--Hood elements.}
\label{fig:Convergence-LDC}
\end{figure}

Figure \ref{fig:Convergence-LDC} shows optimal convergence of $H^1$ velocity error for both classical Taylor--Hood elements of degree $k=2$ and isogeometric Taylor--Hood elements of degree $k=3$.  This supports the hypothesis that the optimal asymptotic convergence proven for the Oseen problem carries over to the steady Navier--Stokes equations.  

\subsubsection{Unsteady Navier--Stokes}\label{sec:2d-tg}

This section studies convergence in an unsteady benchmark, namely the 2D Taylor--Green vortex, using both the dynamic and quasi-static subscale models.  The Taylor--Green vortex has the exact velocity solution
\begin{equation}
    \mathbf{u}(\mathbf{x},t) = \left(\sin(x_1)\cos(x_2)\mathbf{e}_1 - \cos(x_1)\sin(x_2)\mathbf{e}_2\right)e^{-2\nu t}
\end{equation}
on the domain $[-\pi,\pi]^2\times[0,T]$ in space--time, corresponding to a homogeneous source term (i.e., $\mathbf{f}(\mathbf{x},t) = \mathbf{0}$).  This corresponds to free-slip conditions on the spatial boundaries.  In the quasi-static subscale model, we extend the smoothed definition of $\tau_\text{M}$ from \eqref{tau_M_LDC} to
\begin{equation}
\tau_\text{M}=\left(\frac{4}{{\Delta t}^2}+\textbf{u}^h\cdot\textbf{G}\cdot\textbf{u}^h+C_{inv}^2\nu^2\textbf{G}:\textbf{G}\right)^{-\frac{1}{2}}\text{ ,}
\end{equation}
where $\Delta t$ is the time step.

\begin{figure}[!ht]\centering
\includegraphics[width=0.85\textwidth]{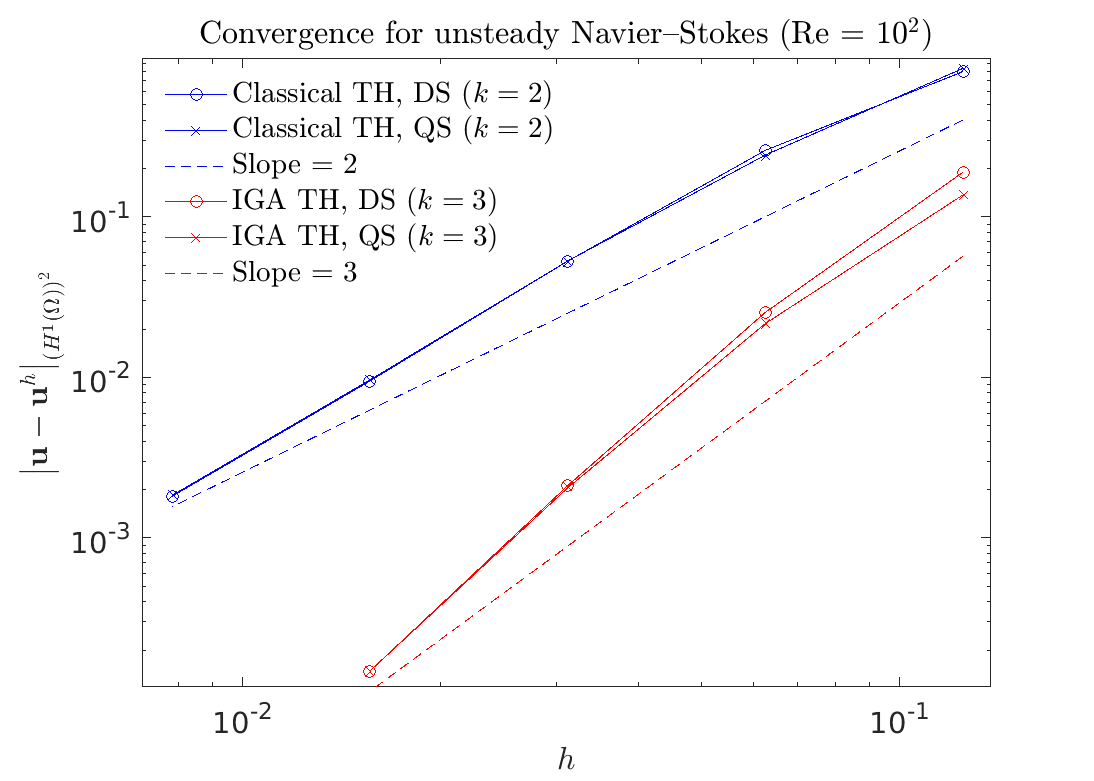}
\caption{$H^1$ velocity error for the 2D Taylor--Green vortex, using classical and isogeometric Taylor--Hood elements, with dynamic subscales (DS) and quasi-static subscales (QS).}
\label{fig:ns-conv-tg}
\end{figure}

Figure \ref{fig:ns-conv-tg} shows the convergence of the $H^1$ seminorm of velocity error at time $T=1$, using dynamic and quasi-static subscale models, in combination with classical Taylor--Hood elements of degree $k=2$ and isogeometric Taylor--Hood elements of degree $k=3$.  This comparison is made at a Reynolds number of 100, in the viscosity-dominated regime, where a direct extrapolation of Theorem \ref{thm:conv} suggests we should expect optimal convergence rates.  The time integration is the implicit midpoint rule with $h \sim \Delta t$ in both cases, which should, in principle, limit convergence rates to second-order, but error is dominated by spatial discretization, and optimal convergence rates are exceeded in both cases.  

\section{Conclusions}\label{sec:conclusions}
We have affirmatively answered the main unresolved conjecture of \cite{Evans2019a}, proving that its proposed stabilized formulation is indeed convergent (in a P\'{e}clet-number-robust way) under the more general conditions that velocity and pressure spaces satisfy a standard inf-sup stability condition.  In \cite[Section 4.3]{Evans2019a}, it was shown numerically that the proposed formulation had some potential for use as a turbulence model when applied with divergence-conforming B-spline discretizations.  We leave to future work the study of whether similar empirical results can be obtained using finite element or isogeometric discretizations that are merely inf-sup stable, such as those used in the convergence tests of this paper.  In extending our discretization of unsteady Navier--Stokes to high Reynolds number flows, it is also worth exploring alternate formulations of advection, such as the ``EMAC'' form \cite{Charnyi2019,Olshanskii2020}, which can retain energetic stability with superior accuracy to the skew form employed in this work.  However, a careful comparison of advective forms is beyond the scope of the present study.  

\section*{Acknowledgements}
SLC worked on this project as part of an independent study course, during a Master's program at UC San Diego.  DK was partially supported by NSF award number 2103939 and JAE was partially supported by NSF award number 2104106.

\bibliographystyle{unsrt}
\bibliography{main}

\end{document}